\documentclass{amsart}
\usepackage{amscd, amsmath, amsthm, amssymb}

\newtheorem{theorem}{Theorem}[section]
\newtheorem{lemma}[theorem]{Lemma}
\newtheorem{proposition}[theorem]{Proposition}

\theoremstyle{definition}     
\newtheorem{definition}[theorem]{Definition}

\newtheorem{question}[theorem]{Question}

\theoremstyle{remark}
\newtheorem{remark}[theorem]{Remark}

\numberwithin{equation}{section}

\begin{document}

\title[Non-commutative free group]
{Tits alternative in hypek\"ahler manifolds}

\author[K. Oguiso]{Keiji Oguiso}
\address{Graduate School of Mathematical Sciences,
University of Tokyo, Komaba, Meguro-ku,
Tokyo 153-8914, Japan, and Korea Institute for Advanced Study, 
207-43 Cheonryangni-2dong, Dongdaemun-gu, Seoul 130-722, Korea}
\email{oguiso@ms.u-tokyo.ac.jp}

\dedicatory{Dedicated to Professor Yukihiko Namikawa on the occasion of
his sixtieth birthday}


\begin{abstract} We show an analogous result of the famous Tits alternative for a group $G$ of birational 
automorphisms 
of a projective hyperk\"aher manifold: Either $G$ contains a non-commutative 
free group or $G$ is an almost abelian group of finite rank. As an application, we show that 
the automorphism groups of the so-called singular K3 surfaces contain non-commutative 
free groups. 
\end{abstract}

\maketitle


\setcounter{section}{0}
\section{Introduction - Background and main results} 
Our main results are Theorems (1.1)-(1.3) and (2.6). Though our actual proof is entirely algebraic, 
these results are much motivated by recent important works about automorphisms of K3 surfaces 
from the view of complex dynamics, notably, a work of McMullen [Mc] and works of Cantat [Ca1, 2].
\par \vskip 1pc  
\noindent 
{\bf 1.} Throughout this note, we work over the complex number field 
$\mathbf C$. By a {\it hyperk\"ahler manifold} 
we mean a compact simply-connected complex K\"ahler manifold 
$M$ having everywhere non-degenerate holomorphic $2$-form $\sigma_{M}$ 
s.t. $H^{2}(M, \Omega_{M}) = \mathbf C \sigma_{M}$. Such manifolds are even dimensional 
and form one of the three fundamental building blocks of compact K\"ahler manifolds 
of vanishing first Chern class. K3 surfaces are nothing but hyperk\"ahler 
manifolds of dimension $2$. We also note that the Hilbert scheme ${\rm Hilb}^{n}\,(S)$,  
of $0$-dimensional closed subschemes of length $n$ of a K3 surface $S$, 
is a hyperk\"ahler manifold of dimension $2n$. We refer to the readers the excellent account 
[JHG, Part III] by Huybrechts about basics of hyperk\"aher manifolds. 
All of what we need here are also reviewed in [Og1, Section 2]. 
We denote the bimeromorphic (resp. biholomorphic) 
automorphism group of 
$M$ by ${\rm Bir}\, M$ (resp. ${\rm Aut}\, (M)$). By an {\it almost abelian group of finite rank}, say $r$, 
we mean a group 
$G$ which containing a subgroup $A$ of finite index which fits in with the exact sequence 
$$1 \longrightarrow F \longrightarrow A \longrightarrow \mathbf Z^{r} \longrightarrow 0\,\, .$$ 
(See eg. [Og1, Section 9].) The extreme counter part of an almost abelian group is a group which contains 
a non-commutative free group, 
or equivalently, a group having a subgroup isomorphic to the free product $\mathbf Z * \mathbf Z$. 
\par \vskip 1pc  
\noindent 
{\bf 2.} 
The aim of this note is to study groups of birational automorphisms of a projective hyperk\"ahler manifold. 
Note that any bimeromorphic automorphism is birational when a manifold is 
projective. Our main results 
are as follows\footnote{The second Theorem (1.2) is inspired by a discussion with Y. Kawamata 
about remarkable difference of the appearance of crescent in Japan and Korea; it looks vertical in Japan 
while horizontal in Korea. This reminded me of the situation in Theorem (1.2).} 
: 
\begin{theorem} \label{theorem:main} 
Let $M$ be a projective hyperk\"ahler manifold and $G$ be a subgroup of ${\rm Bir}\, (M)$. 
Then $G$ satisfies either:
\begin{list}{}{
\setlength{\leftmargin}{10pt}
\setlength{\labelwidth}{6pt}
}
\item[(1)] $G$ is an almost abelian group of finite rank, or
\item[(2)] $G$ contains a non-commutative free group. 
\end{list} 
In particular, so are ${\rm Bir}\,(M)$ and ${\rm Aut}\, (M)$. 
\end{theorem} 
\begin{theorem} \label{theorem:moon} {\rm (Crescent Theorem)}
Let $M$ be a hyperk\"ahler manifold having at least two complex torus fibrations, 
say $\varphi_{i} : M \longrightarrow B_{i}$ ($i = 1$, $2$).  
Then $M$ is projective. Assume furthermore that both $\varphi_{i}$ are Jacobian fibrations of 
positive Mordell-Weil rank. 
Then ${\rm Bir}\,(M)$ contains a non-commutative free group. 
\end{theorem} 
Here by a complex torus fibration, we mean a surjective morphism $f : M \longrightarrow B$ 
over a normal projective variety $B$, whose general fiber is a positive dimensional complex torus. 
(See [Ma] for fibered hyperk\"ahler 
manifolds.) We call $f$ a Jacobian fibration if $M$ is projective and $f$ admits a rational section, say $O$. 
The rational sections of a Jacobian fibration $f$ naturally form an abelian subgroup of 
${\rm Bir}\, (M)$. We call this group the {\it Mordell-Weil group} of $f$ and denote it by 
${\rm MW}\,(f)$. By abuse of language, we call ${\rm MW}\,(f)$ {\it of positive rank} if it contains 
an element of infinite order. See also [Ca1, 2], [CF], [Og2] for relevant 
results. 
\begin{theorem} \label{theorem:sing} 
Let $S$ be a singular K3 surface. Put $M := {\rm Hilb}^{n}\,(S)$. 
Then, both ${\rm Aut}\, (S)$ and ${\rm Aut}\,(M)$, whence ${\rm Bir}\, (M)$, 
contain a 
non-commutative free group. 
\end{theorem} 
Here a {\it singular} K3 surface means a K3 surface having maximum Picard number 
$20$ ([SI]). Theorem (1.3) is a generalization of a famous result of Shioda and Inose 
[ibid], that 
{\it every singular K3 surface has an automorphism of infinite order}, in a non-commutative direction. 
Note that there is a projective K3 surface $S$ s.t. $\rho(S) = \rho$ 
and $\vert {\rm Aut}\, (S) \vert < \infty$ for each integer $\rho$ with 
$1 \le \rho \le 19$ ([Ko, Ni]). So 
a statement similar to Theorem (1.3) is no more true for K3 surfaces of 
smaller Picard number $\rho \le 19$. Theorem (1.3) also shows that the second alternative 
in Theorem (1.1) really occurs in each dimension. 
\par \vskip 1pc  
\noindent 
{\bf 3.} It is also interesting to compare these results with the following 
existing result and interesting open Questions:
\begin{theorem} \label{theorem:non-proj}{\rm [Og1]} 
Let $M$ be a non-projective kyperk\"ahler manifold. Then both ${\rm Bir}\, (M)$ 
and ${\rm Aut}\, (M)$ are almost abelian groups of finite rank. In particular, 
they are finitely generated.
\end{theorem} 
\begin{question} \label{question:fingen} 
Let $M$ be a projective kyperk\"ahler manifold. Are ${\rm Bir}\, (M)$ 
and ${\rm Aut}\, (M)$ finitely generated?
\end{question} 
\begin{question} \label{question:hwang} {\rm (J.M. Hwang)} 
Let $M$ be a projective kyperk\"ahler manifold. Is ${\rm Aut}\, (M)$ 
is of finite index in ${\rm Bir}\, (M)$?
\end{question} 
In dimension $2$, the first question is affirmative by Sterk [St] 
and the second one is trivial. The second question 
is asked by J.M. Hwang after my talk on a relevant subject at KIAS on 
March 2005. 
\par \vskip 1pc  
\noindent 
{\bf 4.} Our proof of Theorems (1.1)-(1.3) is based on two famous, deep results in linear 
algebraic groups; 
Lie-Kolchin Theorem and Tits Theorem ([Hm], [Ti]; see also [Ha] and Section 2 
for the statements) and the notion of Salem polynomial (see eg. [Mc] and also Section 2). Unfortunately, 
our proof does not tell us much about 
algebro-geometrical reason why non-abelian free groups should be in the 
birational automorphism groups in Theorem (1.2). 
It would be interesting to find a more "geometrically visible" proof of Theorem (1.2), especially 
when $M$ is a K3 surfrace. For this, an observation of Cantat [Ca2] might give us some hint. 
\par \vskip 1pc  
\noindent   
{\it Acknowledgement.} An initial idea of this note has been grown up during my stay 
at KIAS March 2005. I would like to express my thanks to Professor 
Y. Kawamata for his 
valuable suggestions and to Professors T.C. Dinh, J.M. Hwang, J.H. Keum, 
S. Morita and D.-Q. Zhang for their interest in this work. I would like to 
express my thanks to Professors J.M. Hwang and B. Kim for invitation. 
Last but not least at all, 
I would like to express my deep thanks to Professor S. Cantat for his several 
valuable comments 
on an earlier version of this note. 

\section{Algebraic preparations} 

The goal of this section is Theorem (2.6), the technical heart in the proof of our main results. 
In {\bf 1}, we recall Lie-Kolchin Theorem and Tits Theorem (Tits alternative). Both are very important 
for our proof. We recall basic notions about lattices in {\bf 2} 
and a few facts about Salem polynimials in {\bf 3}. We show Theorem (2.6) in {\bf 4}. 
\par \vskip 1pc  
\noindent 
{\bf 1.} For simplicity, we shall work only over $\mathbf C$. 
Let $V \not= \{0\}$ be a finite dimensional vector space over $\mathbf C$. 
We regard the general linear 
group ${\rm GL}\,(V)$ as an algebraic group defined over $\mathbf C$, 
with Zariski topology. We identify ${\rm GL}\,(V)$ 
with the group ${\rm GL}\,(V)(\mathbf C)$ of $\mathbf C$-valued points 
in a usual way. 
A subgroup of ${\rm GL}\, (V) = {\rm GL}\, (V)(\mathbf C)$ simply means a subgroup 
as an abstract group. If $G$ is a subgroup of ${\rm GL}\, (V)$, 
then its Zariski closure $\overline{G}$, as well as the identity component of $\overline{G}$, is 
an algebraic subgroup of 
${\rm GL}\, (V)$. First,
we notice the following well-known:
\begin{lemma} \label{lemma:sol}  
Let $G$ be a solvable subgroup of ${\rm GL}\,(V)$. Then: 
\begin{list}{}{
\setlength{\leftmargin}{10pt}
\setlength{\labelwidth}{6pt}
}
\item[(1)] Any subgroup of $G$ and any quotient group of $G$ are solvable.
\item[(2)] $\overline{G}$ is also solvable. 
\end{list} 
\end{lemma}
\begin{proof} We only show (2). It suffices to check that 
$\overline{[G, G]} = [\overline{G}, \overline{G}]$. Note that 
$[\overline{G}, \overline{G}]$ is closed in $\overline{G}$ (See for instance 
[Hm, 17.2]). Thus $\overline{[G, G]} \subset [\overline{G}, \overline{G}]$.

Let us show the other inclusion. Take $g \in G$. Let us define the map $\alpha_{g}$ by
$$\alpha_{g} : \overline{G} \longrightarrow \overline{G}\,\, ;\,\, f \mapsto f^{-1}g^{-1}fg\,\, .$$
Clearly, $\alpha_{g}$ is continuous and satisfies $\alpha_{g}(G) 
\subset [G, G]$. 
Thus $\alpha_{g}(\overline{G}) \subset \overline{[G, G]}$. Hence 
$[\overline{G}, G] 
\subset \overline{[G, G]}$. Let $f \in \overline{G}$. Let us define the map $\beta_{f}$ by 
$$\beta_{f} : \overline{G} \longrightarrow \overline{G}\,\, ;\,\, g \mapsto f^{-1}g^{-1}fg\,\, .$$
Clearly, $\beta_{f}$ is continuous and satisfies $\beta_{f}(G) 
\subset [\overline{G}, G]$. Since $[\overline{G}, G] 
\subset \overline{[G, G]}$, we have $\beta_{f}(G) \subset \overline{[G, G]}$ 
as well. Thus $\beta_{f}(\overline{G}) \subset \overline{[G, G]}$ 
and hence $[\overline{G}, \overline{G}] \subset \overline{[G, G]}$. 
\end{proof}
Lie-Kolchin Theorem and Tits alternative are the following: 

\begin{theorem} \label{theorem:liekol} {\rm (Lie-Kolchin Theorem, see eg. 
[Hm, Chap. VII, 17.6])} 
Let $G$ be a connected solvable subgroup of ${\rm GL}\,(V)$. 
Then $G$ has a common eigenvector in $V$. 
\end{theorem} 

\begin{theorem} \label{theorem:tit} {\rm (Tits alternative [Ti])} 
Let $G$ be a subgroup of ${\rm GL}\,(V)$. 
Then $G$ is either virtually solvable or 
contains a non-commutative free group. 
\end{theorem}

Here a group $G$ is called {\it virtually solvable} if $G$ contains a solvable 
subgroup of finite index. We also notice that any non-commutative free 
group contains $\mathbf Z * \mathbf Z$, i.e. 
the free group of rank $2$. 
\par \vskip 1pc  
\noindent 
{\bf 2.} By a {\it lattice} $L = (L, (*, **))$, we mean a pair consisting 
of a free 
abelian group $L \simeq \mathbf Z^{r}$ and its (possibly degenerate) integral-valued 
symmetric bilinear form 
$(*, **) : L \times L \longrightarrow \mathbf Z$. By $L_{K}$, we denote the scalar extension 
$L \otimes_{\mathbf Z} K$ of $L$ by a field $K$. The signature of $L$ is the pair of the numbers of positive-, 
zero- and negative-eigenvalues of a symmetric matrix associated to $(*,**)$. We call $L$ {\it hyperbolic} (resp. 
{\it parabolic}, {\it elliptic}) if the signature is $(1, 0, r-1)$ (resp. $(0, 1, r-1)$, $(0, 0, r)$). 
\par \vskip 1pc  
We call an element $v \in L \setminus \{0\}$ (resp. a sublattice $M$) {\it primitive} 
if $L/\langle v \rangle$ (resp. $L/M$) 
is torsion-free. We denote by ${\rm O}\, (L)$ 
the group of isometries of a lattice $L$. Note that $\vert {\rm O}\,(L) \vert < \infty$ 
if $L$ is elliptic. 
\par \vskip 1pc  
{\it From now on until the end of this section, we choose and fix a hyperbolic lattice $L$ of rank  
$r$}. 
\par \vskip 1pc  
Then the set $\{ v \in L_{\mathbf R} \vert (v^{2}) > 0\}$ 
consists of two connected components (w.r.t. Euclidean topology of 
$L_{\mathbf R}$). We choose 
and fix one of them and denote it by $\mathcal P(L)$. We call 
$\mathcal P(L)$ the {\it positive cone} of $L$. In general, 
there is no canonical choice of the positive cone. (However, when $L$ is the N\'eron-Severi 
group of a projective hyperk\"ahler manifold, with Beauville-Bogomolov-Fujiki's 
bilinear form, we always choose the positive cone so that it contains ample classes.) 
\par \vskip 1pc  
Let $\overline{\mathcal P}(L)$ (resp. $\partial \mathcal P(L)$) be the 
closure (resp. the boundary) 
of the positive cone in $L_{\mathbf R}$ (w.r.t. Euclidean topology). 
By the Schwartz inequality, we have $(x, y) \ge 0$ for $x, y \in \overline{\mathcal P}(L) \setminus \{0\}$ 
and the equality holds iff $\mathbf R_{>0}x = \mathbf R_{>0}y \in \partial \mathcal P(L)$. 
\par \vskip 1pc  
Let $M$ be a primitive sublattice of $L$. Then $M$ is either hyperbolic, parabolic, or elliptic, 
and $M_{L}^{\perp}$ is elliptic, parabolic, hyperbolic respectively. 
Here and here after we denote by $M_{L}^{\perp}$, the 
primitive sublattice $\{\,v \in L\, \vert\, (v, w) = 0\,\, \forall w \in M\,\}$. Note that 
$M \cap M_{L}^{\perp} = \{0\}$ 
and $M \oplus M_{L}^{\perp}$ is of finite index in $L$ when $M$ is hyperbolic or elliptic, 
while $M \cap M_{L}^{\perp} = \mathbf Z e$ with $(e^{2}) = 0$ and $M + M_{L}^{\perp}$ 
is corank $1$ in $L$ when $M$ is parabolic. 
\par \vskip 1pc  
\noindent  
{\bf 3.} In [Mc] and [Og1], the notion of Salem polynomial plays a very crucial role. In our proof, 
it also plays an important role. We recall the definition:
\begin{definition} \label{definition:salem} An irreducible monic polynomial 
$\Phi(t) \in \mathbf Z[t]$ of degree $n$ is called a Salem polynomial if the complex roots of $\Phi(t) = 0$ 
consists of two real roots $a$ and $1/a$ s.t. $a > 1$ and $n-2$ roots on the unit circle 
$S^{1}$. 
\end{definition} 

In our argement, we need the following result. This is a formal generalization of [Mc, Proposition 3.3]:
\begin{proposition} \label{proposition:salempoly} 
Let $L$ be a hyperbolic lattice of rank $r$ and $g$ be an element of ${\rm O}(L)$, preserving a positive 
cone $\mathcal P(L)$. 
Then the irreducible factors of the characteristic polynomial $\Phi_{g}(t) := {\rm det}\, (tI - g)$ of $g$ 
includes at most one Salem polynomial; and the remaining factors are cyclotomic. In particular, if $c$ 
is an eigenvalue of $g$, then so is $1/c$.  
\end{proposition}
\begin{proof} Since $\Phi_{g}(t) \in \mathbf Z[t]$ is monic, all the eigen values of $g$ 
are algebraic integers. Since $L$ is hyperbolic and $g \in {\rm O}\,(L)$, it follows 
that $g$ has at most one eigenvalue (counted with multiplicity) 
{\it outside} $S^{1}$. Therefore by $\vert {\rm det}\, g\vert = 1$, $g$ satisfies either that all eigenvalues 
are on $S^{1}$, or that $g$ has exactly one eigenvalue outside $S^{1}$, say $a$, exactly one eigenvalue 
inside $S^{1}$, say $b$, and $r-2$ eigenvalues on $S^{1}$.

In the first case, the eigenvalues are all roots of unity by Kronecker's Theorem. So, the irreducible factors 
of $\Phi_{g}(t)$ are cyclotomic.

Let us consider the second case. 
Both $a$ and $b$ are real 
by the uniqueness. Then, one can choose an eigenvector $v$ 
with eigenvalue $a$ in $L_{\mathbf R}$. Then, by $(v^{2}) = (g(v)^{2}) = a^{2}(v^{2})$, we have 
$(v^{2}) = 0$. Thus, $v \in \mathcal P(L)$ (after replacing $v$ by $-v$ if necssary). The same holds 
for an eigenvector with eigenvalue $b$. 
Thus $a > 1$ and $ab = 1$ by $g(\mathcal P(L)) = \mathcal P(L)$ and 
$\vert {\rm det}\, (g) \vert = 1$. Thus ${\rm det}\, g = 1$ as well. Let $f(t) \in \mathbf Z[t]$ be the 
minimal monic polynomial 
of $a$. Then $f(t) \vert \Phi_{g}(t)$. Since $\Phi_{g}(0) = \pm 1$, we have also $f(b) = 0$. 
Thus $f(t)$ is a Salem polynomial. The zeros of the monic polynomial $\Phi_{g}(t)/f(t) \in \mathbf Z[t]$ 
are now on $S^{1}$. Thus, again by Kronecker's Theorem, $\Phi_{g}(t)/f(t)$ includes only cyclotomic polynomials 
as its irreducible factors. The last statement is now clear. Indeed, if $c \in S^{1}$, then 
$1/c = \overline{c}$ is a zero of $\Phi_{g}(t)$. If $c \not\in S^{1}$, then $c = a$ or $b$, 
and $1/c = b$ or $a$ is also a zero of $\Phi_{g}(t)$. 
\end{proof}
\par \vskip 1pc  
\noindent 
{\bf 4.} Let us now formulate our key result:
\begin{theorem} \label{theorem:key} 
Let $L$ be a hyperbolic lattice of rank $r$ and $G$ be a subgroup of ${\rm O}(L)$. 
Assume that $G$ is virtually solvable. Then $G$ is almost abelian of finite rank. 
\end{theorem} 
\begin{proof} We shall show Theorem (2.6) by induction on $r$. 
The result is clear if $r = 1$. So, we may assume that 
$r \ge 2$. We shall proceed the proof by dividing into five steps.
\par \vskip 1pc  
\noindent 
{\bf Step 1.} By the assumption, $G$ has a solvable subgroup 
$N$ s.t. 
$[G : N ] < \infty$. Note that $G$ is almost abelian of finite rank 
iff so is $N$ (cf. [Og1, Section 9]). Then by replacing $G$ by $N$, we may assume that $G$ 
is solvable. 
{\it We will do so from now on}. 
\par \vskip 1pc  
\noindent 
{\bf Step 2.} Put $V := L_{\mathbf C}$. We have a natural embedding:
$$G \subset {\rm O}\,(L) \subset {\rm GL}\,(V)\,\, .$$  
Let $\overline{G}$ be the Zariski closure of $G$ in ${\rm GL}\,(V)$ 
and $S$ be the identity component of $\overline{G}$, i.e. 
the irreducible component containing the identity $1$. Since $G$ is solvable, 
so is $\overline{G}$ by Lemma (2.1)(2). Thus, by Lemma (2.1)(1), 
$S$ is a connected solvable subgroup of ${\rm GL}\, (V)$. Since 
$\overline{G}$ is an algebraic subset of a noetherian space ${\rm GL}\, (V)$, 
it has only finitely many irreducible components. Thus, $[\overline{G} : S] < \infty$. 
For the same reason as in Step 1, we may now assume that $G$ is a subgroup 
of a connected solvable group 
by replacing $G$ by its finite index subgroup 
$G \cap S$, 
and that the positive cone $\mathcal P(L)$ is $G$-stable by further replacing $G$ by its index two subgroup 
(if necessary). 
{\it We will do so from now on}.  
\par \vskip 1pc  
\noindent 
{\bf Step 3.} Recall that $G < {\rm GL}(r, \mathbf Z)$. Then, if all elements of $G$ is of finite order, 
then their orders are universally bounded by some constant, say $N$. Then $G$ is a finite group by 
Burnside's Theorem. In particular, $G$ 
is almost abelian of rank $0$ and we are done. So, we may assume that there 
is $g_{0} \in G$ s.t. ${\rm ord}\, g_{0} = \infty$. {\it We will do so from now on}.  
\par \vskip 1pc  
\noindent 
{\bf Step 4.} Since $G$ is in a connected solvable subgroup, say $S$,  
of ${\rm GL}\, (V)$, by applying Lie-Kolchin Theorem (2.2) for $S$, 
we find a common eigenvector of $G$, say $v \in V = L_{\mathbf C}$. 
Set $g(v) = \alpha(g)v$ for $g \in G$. Then $\alpha$ defines a group homomorphism
$$\alpha : G \longrightarrow \mathbf C^{\times}\,\, ;\,\, g \mapsto \alpha(g)\,\, .$$
Let $M$ be the minimal primitive sublattice of $L$ s.t. $v \in M_{\mathbf C}$. 
This $M$ is $G$-stable. Indeed, since $g \in G$ is defined over $\mathbf Z$ 
and $g(v) = \alpha(g)v$, one has $M \cap g(M) \subset M$, $(M \cap g(M))_{\mathbf C} 
= M_{\mathbf C} \cap g(M)_{\mathbf C}$ and $v \in M_{\mathbf C} \cap g(M)_{\mathbf C}$. 
Thus $v \in (M \cap g(M))_{\mathbf C}$, and therefore $g(M) = M$ by the minimality of $M$. 
\par \vskip 1pc  
\noindent
{\bf Step 5.} Note that $M$ is either elliptic, parabolic, or hyperbolic. 
We completes the proof by dividing into these three cases.
\par \vskip 1pc  
\noindent
{\it The case where $M$ is elliptic.} In this case $M_{L}^{\perp}$ 
is hyperbolic. Let 
$K := {\rm Ker}\, (r_{M} : G \longrightarrow {\rm O}\, (M))$. Then 
$K$ is of inite index in $G$ 
and $K \subset {\rm O}\,(M_{L}^{\perp})$. Since 
${\rm rank}\, M_{L}^{\perp} < {\rm rank}\, L$, 
$K$ is almost abelian of finite rank by the induction hypothesis. 
Thus so is $G$. 
\par \vskip 1pc  
\noindent
{\it The case where $M$ is parabolic.} In this case there is a unique primitive element 
$u \in M$ s.t. $(u^{2}) = 0$ and $u \in \partial\mathcal P(L)$. 
By the uniqueness of $u$ and by $G(\mathcal P(L)) = \mathcal P(L)$, 
this $u$ is $G$-stable. 
Since $u \in \partial \mathcal P(L) \cap L \setminus \{0\}$, it follows from [Og2, Proposition 2.9] 
that $G$ is almost abelian of finite rank.
\par \vskip 1pc  
\noindent
{\it The case where $M$ is hyperbolic.} In this case $M_{L}^{\perp}$ is elliptic 
(possibly $0$).

Consider first the case where $M_{L}^{\perp} \not= \{0\}$, i.e. the case where $M \not= L$. 
Let $K := {\rm Ker}\, (r_{M} : G \longrightarrow {\rm O}\, (M_{L}^{\perp}))$. 
Then $K$ is of finite index in $G$ 
and $K \subset {\rm O}\,(M)$. Since ${\rm rank}\, M < {\rm rank}\, L$ by the case assumption, 
$K$ is almost abelian of finite rank by the induction hypothesis. Thus so is $G$. 

It remains to consider the case where $M = L$. In this case, the result follows from the next Lemma (2.7)(4). 
\begin{lemma} \label{lemma:inj} Assume that $M = L$. Then:
\begin{list}{}{
\setlength{\leftmargin}{10pt}
\setlength{\labelwidth}{6pt}
}
\item[(1)] The homomorphism $\alpha : G \longrightarrow \mathbf C^{\times}$ is injective. In particular, 
$G$ is an abelian group. 
\item[(2)] The characteristic polynomial $\Phi_{g_{0}}(t)$ is a Salem polynomial. Here 
$g_{0}$ is an element of $G$ with ${\rm ord}\, g_{0} = \infty$ in Step 3. 
\item[(3)] Let $g \in G$. Then there is $\varphi(t) \in \mathbf Q[t]$ s.t. $g = \varphi(g_{0})$. 
\item[(4)] $G$ is an almost abelian group of finite rank. 
\end{list}

\end{lemma} 
\begin{proof} Let us show (1). If $\alpha(g) = 1$, then $v \in E(g, 1)$. Here 
$E(g, 1) (\subset L_{\mathbf C})$ is the eigenspace 
of eigenvalue $1$ of $g$. Since $g$ is defined over $\mathbf Z$, there is a primitive sublattice 
$E \subset L$ s.t. $E_{\mathbf C} = E(g, 1)$. Thus, by the minimality of $L = M$, one has $E = L$, i.e. 
$g = 1$. Thus $\alpha$ is injective. 
\par \vskip 1pc  
\noindent
Let us show (2). If $\alpha(g_{0})$ is a root of unity, then there is 
a positive integer $m$ s.t. $\alpha(g_{0}^{m}) = 1$. Then $g_{0}^{m} = 1$ by (1), 
a contradiction to ${\rm ord}\, g_{0} = \infty$. Thus $\alpha(g_{0})$ is not a root of unity. 
Thus, by Proposition (2.5), $\Phi_{g_{0}}(t)$ has a Salem polynomial, say $f(t)$, as its irreducible factor. 
If $\Phi_{g_{0}}(t) \not= f(t)$, then one can write $\Phi_{g_{0}}(t) = f(t)h(t)$. Here $h(t)$ is a product of 
cyclotomic 
polynomials by Proposition (2.5). However, since $(f(t), h(t)) = 1$ and $f(\alpha(g_{0})) = 0$, 
the decomposition $\Phi_{g_{0}}(t) = f(t)h(t)$ would lead a non-trivial rational decomposition 
of $L_{\mathbf Q}$, say $L_{\mathbf Q} = L_{1} \oplus L_{2}$, s.t. $v \in (L_{2})_{\mathbf C}$. 
However this contradicts the minimality of $L = M$. Thus 
$\Phi_{g_{0}}(t) = f(t)$ and the results follows.
\par \vskip 1pc  
\noindent
Let us show (3). Since $\Phi_{g_{0}}(t)$ is irreducible by (2), the eigenvalues of $g_{0}$ are mutaully distinct. 
Thus the $\mathbf C$-linear space $W' := \{h \in {\rm End}_{\mathbf C}(L_{\mathbf C}) \vert hg_{0} = g_{0}h\}$ 
is of dimension $r$. Since $g_{0}$ is defined over $\mathbf Q$, one has $W' = W_{\mathbf C}$, 
where $W := \{h \in {\rm End}_{\mathbf Q}(L_{\mathbf Q}) \vert hg_{0} = g_{0}h\}$. 
Note that $g_{0}^{k}$ ($0 \le k \le r-1$) 
are linearly independent over $\mathbf Q$ by the irreducibility of $\Phi_{g_{0}}(t)$. 
Thus $\langle g_{0}^{k} \rangle_{k=0}^{r-1}$ forms a basis of $W$ over $\mathbf Q$. 
Since $G \subset W$ by (1), the result follows.  
\par \vskip 1pc  
\noindent
Let us show (4). Let $F (\subset \mathbf C)$ be the minimal splitting field of $\Phi_{g_{0}}(t)$ over $\mathbf Q$. 
Let $O_{F}$ be the ring of algebraic integers of $F$ and $U_{F}$ be the unit group of $O_{F}$. 
Then $\alpha(G) \subset F$ by (3). Both $\alpha(g)$ and $1/\alpha(g)$ ($g \in G$) are zero of the 
characteristic polynomial 
$\Phi_{g}(t) \in \mathbf Z[t]$ by Proposition (2.5). Thus both $\alpha(g)$ and $1/\alpha(g)$ are algebraic integers 
in $F$. Hence $\alpha(G) \subset U_{F}$. 
By the Dirichlet unit Theorem, $U_{F}$ is a finitely generated abelian group. Thus so is its 
subgroup $\alpha(G)$. Since $\alpha$ is injective, this implies the result. 
\end{proof}
\par \vskip 1pc  
\noindent
This completes the proof of Theorem (2.6). 
\end{proof} 
\section{Proof of Theorems (1.1)-(1.3)} 
In this section, we shall prove Theorems (1.1)-(1.3). Let $M$ be a projective hyperk\"ahler 
manifold. Then the N\'eron-Severi group $NS(M)$ is a hyperbolic lattice 
w.r.t. the Beauville-Bogomolov-Fujiki's form. Let $G$ be a subgroup of ${\rm Bir}\, (M)$. 
Then, there is a natural group homomorphism $r_{NS} : G \longrightarrow {\rm O}\, (NS(M))$. The kernel 
of $r_{NS}$ is a finite group by [Og1, Corollary 2.7]. So, $G$ contains a non-commutative free group iff so does 
$r_{NS}(G)$, and $G$ is almost abelian of finite rank iff so is $r_{NS}(G)$. 
\par \vskip 1pc  
\noindent 
{\it Proof of Theorem (1.1).} Assume that $r_{NS}(G)$ does not contain a non-commutative free subgroup. 
Then, by the Tits alternative (2.3), $r_{NS}(G)$ is virtually solvable. Thus $r_{NS}(G)$ is almost abelian of 
finite rank by Theorem (2.6). Q.E.D. 
\par \vskip 1pc  
\noindent 
{\it Proof of Theorem (1.2).} Since $B_{i}$ are assumed to be projective, one can write 
$\varphi_{i} = \Phi_{\vert \varphi_{i}^{*}H_{i} \vert}$. 
Here $H_{i}$ is a very ample divisor of $B_{i}$ and $\Phi_{\vert \varphi_{i}^{*}H_{i} \vert}$ 
is the morphism associated with the complete linear system $\vert \varphi_{i}^{*}H_{i} \vert$. Put 
$e_{i} := [\varphi_{i}^{*}H_{i}] \in NS(M)$. 
Then $(e_{i}^{2}) = 0$ (cf. [Ma]). Since $B_{i}$ is $\mathbf Q$-factorial and $\rho(B_{i}) = 1$ by Matsushita 
[ibid], it follows that $\mathbf R_{>0}e_{1} \not= \mathbf R_{>0}e_{2}$. (Here we note that all the arguments in [ibid] 
are valid if we assume that the base space is normal and proective, even if $M$ is not assumed to be projective. 
Indeed, one can rewrite his argument [ibid] by using a K\"ahler class on $M$, instead of an ample class). 
Thus $((e_{1} + e_{2})^{2}) > 0$ and hence $NS(M)$ is hyperbolic. Thus $M$ is projective by the fundamental 
result of Huybrechts [Hu].

Choose $f_{i} \in {\rm MW}\,(\varphi_{i})$ s.t. ${\rm ord}\, (f_{i}) = \infty$. 
We naturally regard both $f_{i}$ as elements of ${\rm Bir}\, (M)$. Set $G := \langle f_{1}, f_{2} \rangle$. 
Then $r_{NS}(G) = \langle r_{NS}(f_{1}), r_{NS}(f_{2}) \rangle$ and 
$r_{NS}(f_{i})(e_{i}) = e_{i}$ for each $i = 1$, $2$. 
By Theorem (1.1), it suffices to check that $r_{NS}(G)$ is not almost abelian of finite rank. 
However, the proof of this fact is the same as [Og2, Theorem 1.6(1)] (except a few obvious modifications). Q.E.D. 
\begin{remark} \label{remark:cantat} 
Our proof of Theorem (1.2) is based on Theorem (1.1) which, as we see above, involves a very deep result, 
the Tits alternative. By using the so-called table-tennis Lemma and some elementary results in hyperbolic geometry 
(see eg. [Ha, Sections 1, 3, 4]), one can also find a non-commutative free group in 
$\langle r_{NS}(f_{1}), r_{NS}(f_{2}) 
\rangle$ more directly. This is 
pointed out to us by S. Cantat. 
\end{remark} 
\par \vskip 1pc  
\noindent 
{\it Proof of Theorem (1.3).} By [Og2, Theorem 1.6(1)], every singular K3 surface $S$ admits at least two 
Jacobina fibrations of positive Mordell-Weil rank. Thus ${\rm Aut}\, (S)$ contains a non-commutative 
free group by Theorem (1.2). By the universality of the Hilbert scheme, ${\rm Aut}\, (S)$ naturally acts on 
$M := {\rm Hilb}^{n}\,(S)$. This action is also faithful. Thus ${\rm Aut}\,(M)$, and hence ${\rm Bir}\,(M)$,  
contains a non-commutative free group. Q.E.D.

\end{document}